\documentclass[12pt]{article}

\usepackage{amsmath,amssymb,euscript,BorelDyn}
\newtheorem{theorem}{Theorem}[section]
\newtheorem{lemma}[theorem]{Lemma}
\newtheorem{proposition}[theorem]{Proposition}
\newtheorem{corollary}[theorem]{Corollary}
\newtheorem{definition}[theorem]{Definition}
\newtheorem{remark}[theorem]{Remark}

\newcommand{\coun}{\medskip\noindent \refstepcounter{theorem}
{\bf\arabic{section}.\arabic{theorem}\;\;}}

\renewcommand{\thetheorem}%
{{\bf\thesection.\arabic{theorem}}}

\begin{document}
\title{Smooth automorphisms and path-connectedness in Borel dynamics}
\author{{\bf S.~Bezuglyi}\thanks{Supported in part by  CRDF grant
UM1-2546-KH-03}, \ \ \ \ {\bf K.~Medynets}\\ Institute for Low
Temperature Physics, Kharkov, Ukraine}
\date{}

\maketitle

\newcommand{\e}{\varepsilon}
\newcommand{\N}{{\mathbb N}}
\newcommand{\au}{Aut(X,\mathcal{B})}
\newcommand{\auo}{Aut_0(X,\mathcal{B})}
\newcommand{\M}{\mathcal{M}_1(X)}
\newcommand{\la}{\lambda}
\newcommand{\per}{{\cal P}er}
\newcommand{\ap}{{\cal A}p}
\newcommand{\A}{\mathcal{A}}
\newcommand{\B}{\mathcal{B}}
%\centerline{\it Draft}

\begin{abstract}
Let $Aut(X,\mathcal{B})$ be the group of all Borel automorphisms of a
standard Borel space $(X,\mathcal{B})$. We study topological properties of
$\au$ with respect to the uniform and weak topologies, $\tau$ and $p$,
defined in \cite{BDK1}. It is proved that the class of smooth automorphisms
is dense in $(\au,p)$. Let $Ctbl(X)$ denote the group of Borel automorphisms
with countable support. It is shown that the topological group
$\auo=Aut(X,\mathcal{B})/Ctbl(X)$ is path-connected with respect to the
quotient topology $\tau_0$. It is also proved that $\auo$ has the Rokhlin
property in the quotient topology $p_0$, i.e. the action of $\auo$ on itself
by conjugation is topologically transitive.
\end{abstract}

%%%%%%%%%%%%%%%%%%%%%%%%%%%%%%%%%%%%%%%%%%%%%%%%%%%%%%%%%%%%%%
%
%                   INTRODUCTION
\section{ Introduction}

In the present paper, we continue the study of topological properties of
the group $\au$ of all Borel automorphisms of a standard Borel space
$(X,\mathcal{B})$. We consider two topologies, $\tau$ and $p$, on $\au$
which take their origins in ergodic theory. They were defined and studied
in the context of Borel and Cantor dynamics in
\cite{BDK1,BDK2,BDM,BK1,BK2}. Recall that the topology $\tau$ is defined by
the base of neighborhoods $U(T;\mu_1,\ldots,\mu_n;\varepsilon)=\{S\in
\au\,|\, \mu_i(E(S,T))<\varepsilon, \;i=1,\ldots,n\}$, where
$\mu_1,\ldots,\mu_n$ are Borel probability measures on $X$ and
$E(S,T)=\{x\in X : Tx\neq Sx\}\cup\{x\in X : S^{-1}x\neq T^{-1}x\}$.
Obviously, $\tau$ is a direct analogue of the well known uniform topology
on the group $Aut(X,\mathcal B,\mu)$ of all non-singular automorphisms of a
measure space generated by the metric $d(S,T) =\mu(E(S,T))$. It is
worthwhile to mention that, in fact, $Aut(X,\mathcal B,\mu)$ is formed by
classes of automorphisms coinciding $\mu$-almost everywhere. It allows one
to neglect the behavior of automorphisms on sets of zero measure.
Topological properties of $(Aut(X,\mathcal B,\mu),d)$ were extensively
studied in ergodic theory (see, for example, \cite{AP,BG,CSF,D,H,R}). The
other topology, $p$, on $\au$ is defined by neighborhoods
$W(T;F_1,\ldots,F_n)=\{S\in Aut(X,\mathcal{B})\ |\ SF_i= TF_i,\;
i=1,\ldots,n \}$ where $F_1,\ldots, F_n$ are Borel sets. It was shown in
\cite{BDK1} that $p$ can be treated as an analogue of the weak topology
$d_w$ which has been also widely used in ergodic theory. Note also that, in
the context of Cantor dynamics, $p$ is equivalent to the $\sup$-topology on
the set of all homeomorphisms. Based on this observation, we  call $\tau$
and $p$ the uniform and weak topologies on $\au$ respectively.

Our goal is to find out which of topological properties known in ergodic
theory for $Aut(X,\mathcal B,\mu)$ hold for $\au$ with respect to the
topologies $\tau$ and $p$. For instance, it is important for many
applications to know dense subsets in $\au$ which consist of
\textquotedblleft relatively simple" Borel automorphisms. It was shown in
[BDK1] that the set $\per$ of periodic automorphisms is dense in $(\au,
\tau)$ but non-dense in $(\au,p)$. Therefore one needs to extend $\per$ to
produce a dense subset in $\au$ with respect to $p$. The class $\mathcal{S}$
of smooth Borel automorphisms is a natural extension of periodic
automorphisms. By definition, $T$ is smooth if there exists a Borel subset
in $X$ which intersects every $T$-orbit exactly once. In this paper (see
Section 2), we prove that the set of smooth Borel automorphisms is dense in
$(\au, p)$. This statement has been conjectured by A.~Kechris.

A number of papers in ergodic theory was devoted to the study of
connectedness of $Aut(X,\mathcal B,\mu)$ in the weak and uniform topologies
(see, for example, \cite{BG,D,H,Ke}). In particular, it was proved that
$Aut(X,\mathcal B,\mu)$ was path-connected and even contractible. It is not
hard to show that $\au$ is not path-connected in $\tau$ because there is no
continuous path connecting the identity with an automorphism with countable
support. At first sight, the situation seems to be different in Borel
dynamics. But if we factorize $\au$ (like in ergodic theory) by a closed
normal subgroup consisting of automorphisms whose behavior can be
neglected, we do produce a path-connected quotient group. From this point
of view, it is natural to say that $S,T\in \au$ are equivalent if they are
different on an at most countable set. If $Ctbl(X)$ denotes the group of
all Borel automorphisms with countable support, then $\auo= \au/Ctbl(X)$ is
a Hausdorff topological group with respect to the quotient topologies
$\tau_0$ and $p_0$. Note that the group $Ctbl(X)$ was also considered in
\cite{Sh}. It turns out that such a factorization improves topological
properties $\au$ in both quotient topologies $\tau_0$ and $p_0$. Namely,
the group $(\auo, \tau_0)$ becomes path-connected and $(\auo, p_0)$
possesses the so called {\it Rokhlin property}, i.e. the action of $\auo$
on itself by conjugation is topologically transitive (see \cite{GW,BDM} for
the Rokhlin property in Cantor dynamics). In fact, we prove an even
stronger result by showing that the conjugacy class of every aperiodic
smooth automorphism is dense.

Most definitions and notions used in this paper are mostly taken from the
book \cite{Nad}. We collected in Section 1 definitions and facts which are
used in the paper. When we say that $T$ is an automorphism of $(X,\B)$, we
always mean that $T$ is a Borel automorphism. We will also use the term
\textquotedblleft automorphism" for elements of the quotient group $\auo$.

%%%%%%%%%%%%%%%%%%%%%%%%%%%%%%%%%%%%%%%%%%%%%%%%%%%%%%%%%%%%%%%%%%

\section{Preliminaries}

\coun Let $(X,\mathcal{B})$ be a {\it standard Borel space} with the
$\sigma$-algebra of Borel sets $\mathcal{B}$. This means, by definition,
that $(X,\B)$ is (Borel) isomorphic to a Polish space, i.e. a complete
separable metric space. Recall several facts about standard Borel spaces:
(i) any two standard Borel spaces are Borel isomorphic; (ii) if $A\in \B$,
then $A$ is either at most countable or has the cardinality continuum; (iii)
if Borel sets $A,B$ have the same cardinality, then they are isomorphic.

\coun Denote by $\au$ the group of all Borel automorphisms of
$(X,\mathcal{B})$. Let $T\in \au$ and $A\in \B$. The set
$\bigcup_{n\in\mathbb{Z}}T^nA$ is called the {\it saturation} of $A$ with
respect to $T$ and denoted by $s_TA$ (or simply $sA$ if $T$ is clear from
the context). A Borel set $W$ is said to be {\it wandering} with respect to
$T$ if $T^nW \cap W= \emptyset,\ n\in \N$. A Borel set $A\subset X$ is
called a {\it complete section} with respect to $T$ (or simply a
$T$-section) if every $T$-orbit meets $A$ at least once, i.e $s_TA=X$. A
point $x$ from a Borel set $A$ is called {\it recurrent} with respect to $T$
if there exists $n\in\mathbb{N}$ such that $T^nx\in A$.

\coun Denote by $\ap$ and $\per$ the sets of aperiodic and periodic
automorphisms respectively.

We say that a transformation $T\in Aut(X,\mathcal{B})$ is {\it smooth} if
there exists a complete Borel section $A$ such that $A$ meets every
$T$-orbit exactly once. We will denote the class of smooth automorphisms by
$\mathcal{S}$. Obviously, any periodic Borel automorphism is smooth. On the
other hand, if $X$ is a compact metric space and $T$ is an aperiodic
homeomorphism of $X$, then $T$ cannot be smooth.

\coun\label{technicals} We will use the following basic statements taken
from \cite{Nad}.

\begin{enumerate}
\item[(a)] {\it(Poincar\'e Recurrence Lemma)} Let $T\in Aut(X,\mathcal{B})$
and $A\in\mathcal{B}$. Then there exists a wandering set $W\subset A$ such
that for each $x\in A-\bigcup_{k\in\mathbb{Z}}T^kW$ the point $x$ returns to
$A$ for infinitely many positive $n$ and also for infinitely many negative
$n$.

\item[(b)] Let $T\in Aut(X,\mathcal{B})$. Then $X$ can be partitioned into a
disjoint union of Borel sets $X=X_\infty\cup\bigcup_{k\geq 1}X_k$ where
points from $X_k, \ k<\infty$, have period $k$ and $X_\infty$ consists of
aperiodic points.

\item[(c)] Let $T\in \au$ and let $X_k,\ k< \infty$, be as in (b). Then there
exists a Borel set $B_k\subset X_k$ such that
$X_k=\bigcup_{i=0}^{k-1}T^iB_k$ and the union is disjoint.
\end{enumerate}

\coun {\bf $T$-Towers }\label{Towers_Construction} Let $T\in
Aut(X,\mathcal{B})$. Assume that all points from $A$ are recurrent with
respect to $T$. For $x\in A$, define $n(x)=n_A(x)$ as the smallest positive
integer such that $T^{n(x)}x\in A$ and $T^ix\notin A,\ 0<i<n(x)$. Let
$C_k=\{x\in A\ |\ n_A(x)=k\},\ k\in\mathbb{N}$ (some of the $C_k$'s may be
empty). Notice that $T^kC_k\subset A$ and $\xi_k=
\{T^iC_k\,|\,i=0,\ldots,k-1\}$ consists of pairwise disjoint sets. We call
$\xi_k$ a $T$-{\it tower} with base $C_k$ and top $T^{k-1}C_k$. The set
$T^iC_k$ is called the $i$-th level of $\xi_k$. The height of $\xi_k$ is
$k$.

Since $T^nx\in A$ for infinitely many positive and negative $n$, we have
$$
\bigcup\limits_{n\geq 0}T^nA=\bigcup\limits_{n\in\mathbb{Z}}T^nA= s_TA
$$
and
$$
\bigcup\limits_{n\geq 0}T^nA =
\bigcup\limits_{k=1}^\infty\bigcup\limits_{i=0}^{k-1}T^iC_k.
$$
The above relation shows that $\xi= \{\xi_k : k\in \N\}$ forms a partition of
$s_TA$ into $T$-towers $\xi_k,\ k\in \N$. Notice that $T$ maps the union of
tops of these towers onto the union of their bases.

Given a partition $\xi$ on $X$, a Borel set $B\subset X$ is
called a $\xi$-set if it is a union of atoms of $\xi$.

\coun The next lemma is one of the main tools in the study of Borel
automorphisms. It is used in various problems related to finding a suitable
approximation of an aperiodic automorphism, in particular, in the proof of
the Rokhlin lemma \cite{KM,Nad,BDK1}.

\begin{lemma}\label{markers} Let $T\in \au$ be an aperiodic Borel
automorphism of a standard Borel space $(X,\B)$, i.e. $X=X_\infty$. Then
there exists a sequence $(A_n)$ of Borel sets such that\\ (i) $X=A_0 \supset
A_1 \supset A_2\supset \cdots,$\\ (ii) $\bigcap_n A_n =\emptyset,$\\ (iii)
$A_n$ and $X\setminus A_n$ are complete $T$-sections, $n\in \N$,\\ (iv) for
$n\in \N$, every point in $A_n$ is recurrent,\\ (v) for $n\in \N$,\ $A_n\cap
T^i(A_n) =\emptyset,\ i=1,...,n-1$,\\ (vi) the base $C_k(n)$ of every
non-empty $T$-tower constructed over $A_n$ is an uncountable Borel set,
$n\in \N$.
\end{lemma}

For the proof, see \cite[Lemma 4.5.3]{BeK} where (i) - (iii) have been
proved in more general settings of countable Borel equivalence relations.
It is shown in \cite[Chapter 7]{Nad} that one can refine the  choice of
$(A_n)$ to get (iv) and (v). It is clear that one can remove an at most
countable set of points from each $A_n$ to satisfy property (vi).

A sequence of Borel sets $(A_n)$ satisfying \ref{markers} is called a
{\it vanishing sequence of markers}.

\coun Recall now the definition of the uniform and weak topologies on $\au$
following \cite{BDK1}. Let $\mathcal{M}_1(X)$ denote the set of all Borel
probability measures on $X$. A measure $\mu\in \mathcal{M}_1(X)$ is called
continuous (non-atomic) if $\mu(\{x\}) =0$ for all $x\in X$. The Dirac
measure at $x\in X$ is denoted by $\delta_x$. For $T,S\in
Aut(X,\mathcal{B})$, define $E(S,T)=\{x\in X : Tx\neq Sx\}\cup \{x\in X :
S^{-1}x\neq T^{-1}x\}$.

\begin{definition}\label{top} The topologies $\tau$ and $p$ on
$Aut(X,\mathcal{B})$
are defined, respectively, by the bases of neighborhoods $\mathcal{U}=
\{U(T;\mu_1,\ldots,\mu_n;\varepsilon)\}$ and $\mathcal{W}=
\{W(T;F_1,\ldots,F_n)\}$ where
$$\begin{array}{ll}
U(T;\mu_1,\mu_2,\ldots,\mu_n;\varepsilon)=\{S\in \au\ |\ \mu_i(E(S,T))<
\varepsilon,\;i=1,\ldots,n\},\\
\\
W(T;F_1,F_2,\ldots,F_n)=\{S\in Aut(X,\mathcal{B})|\; SF_i=TF_i,\;
i=1,\ldots,n \}.
\end{array}
$$
Here $T\in\au$, $\mu_1,\ldots,\mu_n\in\M$, $\e>0$, and $F_1,\ldots,F_n\in
\mathcal{B}$.
\end{definition}

It was shown in \cite{BDK1} that $\au$ is a Hausdorff topological group with
respect to these topologies. More topological properties of $\au$ and its
subsets can be found in \cite{BDK1}.

\begin{remark}\label{E0} {\rm If in the definition of $\tau$ one takes the
set $E_0(T,S) = \{x \in X : Sx \neq Tx\}$, then the obtained new topology
is, in fact, equivalent to $\tau$. The proof of this fact is straightforward}.
\end{remark}

 \coun Let $Ctbl(X)$ be
defined as a subset of $\au$ consisting of all automorphisms with countable
support, that is $T\in Ctbl(X)$ if $|\{x\in X : Tx\neq x\}|\leq \aleph_0$
where $|A|$ denotes the cardinality of $A$. Note that $Ctbl(X)$ is a normal
subgroup closed with respect to the topologies $\tau$ and $p$, see
 the  proposition below. Therefore
$\auo=\au/Ctbl(X)$ is a Hausdorff topological group with respect to the
quotient topologies $\tau_0$ and $p_0$. Considering elements from $\auo$, we
will identify Borel automorphisms which are different on a countable set,
that is $S\sim S'$ if $|E(S,S')| \leq \aleph_0$. In other words, $S\sim S'$
if there exists $P\in Ctbl(X)$ such that $S= S'P$. This identification
corresponds to the well known approach used in measurable dynamics when two
automorphisms are also identified if they are different on a set of measure
0.

\begin{proposition}  $Ctbl(X)$ is a normal closed subgroup in $\au$
with respect to the topologies $\tau$ and $p$.
\end{proposition}
{\it Proof.} It is obvious that $Ctbl(X)$ is a normal subgroup in $\au$, so
it is enough to prove that it is closed in $\tau$ and $p$. To do this,
suppose that there exists an automorphism $S\in
\overline{Ctbl(X)}^\tau\setminus Ctbl(X)$. Then for any neighborhood
$U(S)=U(S;\mu_1,\ldots,\mu_n;\e)$ there exists an automorphism $R\in
U(S)\cap Ctbl(X)$, that is $\mu_i(E(R,S))<\e$ for all $i$. Since $S\notin
Ctbl(X)$, we have that $E=\{x\in X : Sx\neq x\}$ is uncountable. Let $\nu$
be a continuous measure on $X$ such that $\nu(X \setminus E)=0$. Consider a
neighborhood $U_1 = U(S;\nu;\e)$ of $S$. Then for $R\in U_1\cap Ctbl(X)$ we
have that
$$
\nu(\{x\in X : Sx\neq x\})=1,\qquad \nu(\{x\in X : Sx\neq Rx\})<\e
$$
But $\nu(\{x\in X : Rx\neq x\})=0$, therefore $\nu(\{x\in X : Sx=x\})=0$ and
$\nu(\{x\in X : Sx=Rx=x\})>1-\e$, a contradiction.

Assume now that $S\in\overline{Ctbl(X)}^p \setminus Ctbl(X)$. Then $E=\{x\in
X : Sx\neq x\}$ is uncountable and $S$-invariant. Let $E_1$ be an
uncountable Borel subset of $E$ such that $SE_1\cap E_1= \emptyset$. Then
for every $R\in W(S; E_1)$ we have that $RE_1\cap E_1=\emptyset$. In
particular, if $R\in W(S; E_1)\cap Ctbl(X)$, then we obtain that $R$ acts
non-trivially on the uncountable set $E_1$. This contradicts the definition
of $Ctbl(X)$. \hfill$\square$
\\

The bases of neighborhoods for $\tau_0$ and $p_0$ consists of the sets
$U_0(T;\mu_1,...,\mu_n;\e) = U(T; \mu_1,...,\mu_n;\e)Ctbl(X)$ and
$W_0(T;F_1,...,F_m) = W(T;F_1,...,F_m)Ctbl(X)$, respectively. The next
proposition shows that, in fact, $\tau_0$ and $p_0$ are generated by
neighborhoods $U_0$ and $W_0$ with continuous measures $\mu_i$ and
uncountable sets $F_j$.

\begin{proposition}\label{coun} Given a $\tau_0$-neighborhood  $U_0=
U_0(T;\mu_1,...,\mu_n;\e)$ and a $p_0$-neighborhood ${W_0=
W_0(T;F_1,...,F_m)}$, there exist neighborhoods $U'_0(T;\nu_1,...,\nu_n;\e)
= U'(T;\nu_1,...,\nu_n;\e)Ctbl(X)$ and $W'_0(T;B_1,...,B_m) =
W'(T;B_1,...,B_m)Ctbl(X)$ of $\tau_0$ and $p_0$, respectively, such that
$U'_0 \subset U_0,\ W'_0\subset W_0$ where measures $\nu_1,...,\nu_n$ are
continuous and Borel sets $B_1,...,B_m$ are uncountable.
\end{proposition}
{\it Proof}. Consider the countable set $A=\bigcup_{i=1}^n A_i$ where $A_i
= \{x\in X : \mu_i(\{x\}) > 0\}$. Let $c_i = \mu_i(A)$ and assume that $c_i
< 1$, $i=1,...,n$. Define
$$
\nu_i(B) = \frac{\mu_i(B\cap A^c)}{\mu_i(A^c)},\ \ B\in \B,\ \ i=1,...,n,
$$
where $A^c := X \setminus A$. Clearly, $\nu_i$ is a non-atomic Borel
probability measure on $X$. It remains to show that $U'_0 =
U'(T;\nu_1,...,\nu_n;\e)Ctbl(X)$ is a subset of $U_0$. To do this, it
suffices to check that for every $S\in U'(T;\nu_1,...,\nu_n;\e)$ there
exists $S_1\in U(T;\mu_1,...,\mu_n;\e)$ such that $S\sim S_1$. Let $\Gamma$
be the countable group of automorphisms of $X$ generated by $T$ and $S$. Let
$D = s_\Gamma A$ be the $\Gamma$-orbit of $A$. Define
$$
S_1x=\left\{ \begin{array}{ll}
Tx, &\ \ \ x\in D\\
Sx, &\ \ \ x\in D^c
\end{array} \right.
$$
Obviously, $S_1\sim S$. Since $E(T,S_1) \subset D^c \subset A^c$ and
$E(T,S_1)\subset E(T,S)$, we have that for $i=1,...,n$,
$$
\mu_i(E(T,S_1)) = \mu_i(A^c)\nu_i(E(T,S_1)) \leq
\mu_i(A^c)\nu_i(E(T,S)) < \e\mu_i(A^c) < \e.
$$

Notice that if $\mu_i(A) =1$ for some $i$, then for any $S\in \au$
there exists $S_1\sim S$ such that $\mu_i(E(S_1,T)) = 0$.

The proof for the topology $p_0$ is similar.\hfill$\square$

%%%%%%%%%%%%%%%%%%%%%%%%%%%%%%%%%%%%%%%%%%%%%%%%%%%%%%%%%%

\section{ Smooth automorphisms are dense in $(\au, p)$}

 \coun In this section, we prove that the $p$-closure of the set
$\mathcal S$ of smooth automorphisms is the entire group $\au$. Moreover, it
is shown that $(\auo,p_0)$ has the Rokhlin property.

\begin{theorem}\label{smooth}
$\overline{\mathcal{S}}^p=Aut(X,\mathcal{B})$. Moreover, each
$p$-neighborhood of an aperiodic automorphism  necessarily contains an
aperiodic smooth automorphism.
\end{theorem}
{\it Proof.} Let $T\in Aut(X,\mathcal{B})$. Obviously, it suffices to
consider the case when $T$ is aperiodic. Take a $p$-neighborhood
$W=W(T;F_1,F_2\ldots, F_n)$ of $T$. Without loss of generality, we can
assume that the sets $\{F_1,F_2\ldots,F_n\}$ form a partition of $X$. Show
that there exists a smooth aperiodic automorphism $S\in W$.

For  $i=1,2,\ldots,n$, consider the Borel sets
$$
F_i^j=F_i\cap TF_j,\qquad j=1,2,\ldots,n.
$$
Suppose that the collection $\{F_i^j\}_{i,j=1}^n$ contains exactly $q\leq
n^2$ non-empty sets. Denote them by $V_l$, $l=1,\ldots,q$. Then
$\{V_1,\ldots, V_q\}$ is a partition of $X$ which refines $\{F_1,\ldots,
F_n\}$.

Let
$$
\begin{array}{lll}
B_1=V_1\\
B_2=V_2-sB_1\\
\cdots\cdots\cdots\cdots\\
B_q=V_q-s(B_1\cup B_2\cup\ldots\cup B_{q-1}).
\end{array}
$$
Here and below $s$ stands for $s_T$. Clearly, $\{sB_1,sB_2,\ldots,sB_q\}$ is
a partition of $X$ into Borel sets. Without loss of generality, we may
assume that all sets $B_1,...,B_q$ are non-empty. For each $B_i$, find a
wandering subset $A_i\subset B_i$ (see \ref{technicals}) such that all
points from $D_i= B_i - sA_i$ are recurrent. Therefore, by
\ref{Towers_Construction}, we can find a partition $\Xi_i$ of $sD_i$ into
pairwise disjoint $T$-towers $\{\xi_i(k) : k\in \mathbb N\}$ such that the
union of bases of these towers is $D_i$. For short, we will write
$$
sD_i=\bigcup\limits_{\xi\in\ \Xi_i}\xi,\ \ i=1,...,q.
$$
For a tower $\xi$, denote by $B_\xi$ and $h_\xi$ its base and height,
respectively. By construction, we have that
$$
\bigcup\limits_{\xi\in \Xi_i}B_\xi= D_i \subset V_i\subset F_m
$$
and
$$\bigcup\limits_{\xi\in\Xi_i}T^{h_\xi-1}B_\xi= T^{-1}D_i\subset
T^{-1}V_i=T^{-1}(F_k\cap TF_l)\subset F_{l}
$$
for some $1\leq m, l\leq n$.

Since each $sB_i,\ i=1,...,q$, can be represented as a disjoint union of
$sA_i$ and $sD_i$, we obtain the partition
$$
X=\bigcup\limits_{i=1}^q sB_i=\bigcup\limits_{i=1}^q sA_i\cup
\bigcup\limits_{i=1}^q sD_i.
$$
The automorphism $T$ restricted to $s\bigcup_{s=1}^q A_i=\bigcup_{i=1}^q
sA_i$ is smooth since $\bigcup_{i=1}^qA_i$ is a wandering set for $T$.
Define $S=T$ on $s\bigcup_{s=1}^q A_i$. To complete the proof, we need to
define $S$ on $\bigcup_{i=1}^q sD_i$ such that $S(sD_i\cap F_j)=T(sD_i\cap
F_j)$ for all $j=1,\ldots,n,\ i=1,...,q$.

Fix the set $sD_i$ and let $\Xi_i^0=\{\xi\in \Xi_i : |B_\xi| \leq\aleph_0\}$.
Then
$$
X_i^0=\bigcup\limits_{\xi\in\Xi^{0}_i}sB_\xi
$$
is countable and therefore $T$, restricted to $X^0_i$, is smooth. Set $S =
T$ on $X_0^i$.

For $\xi\in \Xi_i$, denote by $\xi':= \{T^j (B_\xi-X_0^i) : 0\leq j\leq
h_\xi-1\}$. Setting $\Xi'_i=\{\xi' : \xi\in \Xi_i\}$, we get a disjoint
union
$$
sD_i=X_0^i\cup\bigcup\limits_{\xi'\in\Xi'_i}\xi'.
$$
Note that the cardinality of each tower $\xi'\in\Xi'_i$ is continuum.

Now, we define the automorphism $S$ on each tower $\xi'\in \Xi'_i$ such that
$S(\xi') = \xi'$ and $S$ coincides with $T$ on each level of the tower
$\xi'$ except the top of $\xi'$. To do this, we write down the base
$B_{\xi'}$ of $\xi'$ as a disjoint union $B_{\xi'}=\bigcup_{m\in\mathbb{Z}}
B_{\xi'}(m)$ with uncountable Borel sets $B_{\xi'}(m)$. Let $R_m$ be an
arbitrary Borel isomorphism between $T^{h_{\xi'}-1}B_{\xi'}(m)$ and
$B_{\xi'}(m+1),\ m\in \mathbb Z$ (one can assume that $h_{\xi'}\geq 2$ by
\cite[Theorem 7.25]{Nad}). Define $Sx = Tx$ for $x\in \{T^jB_{\xi'}(m) :
0\leq j\leq h_\xi-2,\,m\in \mathbb{Z}\}$ and $Sx = R_mx$ for $x\in
T^{h_\xi-1}B_\xi(m),\ m\in \mathbb Z$. Then $S$ is defined everywhere
on $X$. To prove that $S\in W$, we need to show that $SF_j= TF_j$.
Write down $F_j\cap sD_i$ as a disjoint union of sets $E_0$ and $E_1$
where $E_1=\{x\in F_j\cap sD_i : x\in T^{h_\xi-1}B_\xi,\mbox{ for some }
\xi\in \Xi_i\}$ and $E_0=(F_j\cap sD_i)-E_1$. It follows from the
construction that if $E_1\neq\emptyset$,
then $E_1\subset F_j$ and $E_1=\bigcup_{\xi\in \Xi_i} T^{h_\xi-1}B_\xi =
T^{-1}D_i$. Therefore, by definition of $S$, we have that $SE_1=TE_1$. It is
clear that $S= T$ on $E_0$. The proof is complete. \hfill$\square$

\begin{remark}\label{den} {\rm (1) We note that the set $\mathcal{S}$ is
dense in $(\au,\tau)$. It follows from the fact that the set of periodic
automorphisms is dense in $(\au,\tau)$ \cite[Corollary 2.6]{BDK1}. On
the other hand, $\mathcal{S}\cap\ap$ is nowhere dense in $(\au,\tau)$ by
\cite[Theorem 2.8]{BDK1}.

(2) The set of aperiodic smooth automorphisms is not  dense in $(\au, p)$
because $\ap$ is a closed subset in $(\au,p)$ \cite[Theorem 2.8]{BDK1}}.
\end{remark}

\begin{remark} {\rm B.~Miller proved the following result \cite{M} which
may be used to simplify the proof of Theorem \ref{smooth}:

{\it Suppose that $X$ is a Polish space, $T: X\to X$ is a Borel
automorphism, and $\{A_n\}_{n\in \N}$ is a partition of $X$ into Borel sets.
Then there is a countable $T$-invariant set $C\subset X$ and an aperiodic
smooth Borel automorphism $S : X\setminus C \to X\setminus C$ such that
$T(A_n\setminus C) = S(A_n \setminus C)$ for every $n\in \N$.}}
\end{remark}

\begin{proposition}\label{nowd}  $\per$  is a closed
nowhere dense subset of $\au$ with respect to $p$.
\end{proposition}
{\it Proof}. We first show that the set $\per$ is closed. If we assume that
there exists an automorphism $T\in\overline{\per}^p \setminus\per$, then $T$
must have an aperiodic point $x_0\in X$. Define $F=\{x_0\}\cup\{Tx_0\}\cup
\{T^2x_0\} \cup\ldots$ and consider the $p$-neighborhood $W=W(T;F)$. Then
$W$ necessarily contains a periodic automorphism $P$ such that $PF=TF$.
Since $TF\subsetneqq F$, we obtain that $P^nF\subsetneqq F$ for all
$n\in\N$. Therefore the point $\{x_0\}=F \setminus PF$ must be aperiodic for
$P$, a contradiction.

Let $W(P)=W(P;F_1,\ldots,F_n)$ be a $p$-neighborhood of a periodic
automorphism $P$. We can assume that the sets $(F_1,...,F_n)$ form a
partition of $X$. We will first show that $W(P)$ contains a non-periodic
automorphism. By \ref{technicals}, $X$ is partitioned into $P$-towers
$\xi_k =\{B_k,\ldots,P^{k-1}B_k\}$ such that $P$ has period $k$ on the set
$X_k = B_k\cup \cdots \cup P^{k-1}B_k$. One can refine the partition $\xi =
(\xi_k : k\in \N)$ to produce a new partition $\xi'$ such that every $F_i,
i=1,...,n,$ is a $\xi'$-set. Let $(B,...,T^{m-1}B)$ be a $P$-tower from
$\xi'$ with uncountable base. As in \ref{smooth}, we can find an aperiodic
automorphism $T$ defined on $C =\bigcup_{i=0}^{m-1} T^iB$ such that
$T(P^iB) = P^{i+1}B,\ i=0,...,n-2,$ and $T(P^{n-1}B) = B$. Define $T$ on
$X\setminus C$ by setting $T = P$. We see that $TF_i = PF_i$ for all $i$,
i.e. $T\in W(P)$. It is clear that there exists a Borel set $F\subset C$
such that $TF \subsetneqq F$. Then $W_1 = W(T; F)$ contains no periodic
automorphism. Thus, $(W(T) \cap W(P)) \cap \per = \emptyset$ and we are
done. \hfill$\square$

\coun In contrast to \ref{den}, the situation for the quotient group $\auo$
 is different. It turns out that the set $\mathcal{S}\cap\ap$ is a dense
subset in $(\auo,p_0)$.

\begin{theorem}\label{dense} The set of aperiodic smooth automorphisms
is dense in $(\auo,p_0)$.
\end{theorem}
{\it Proof.} By  \ref{smooth}, we only need to show that each
neighborhood $W_0=W_0(P;F_1,\ldots,F_n)$ of $P\in\per$ contains an aperiodic
smooth automorphism $S$. It follows from \ref{technicals} that $X$ is
decomposed into an at most countable collection of $P$-invariant towers $\Xi
=\{\xi_k : k\in\N\}$. By  \ref{coun}, we can assume that all
$\xi_k$'s are uncountable Borel sets. Let $\Xi'$ be a refinement of $\Xi$,
obtained by cutting the towers from $\Xi$, such that each $F_i$ is a $\Xi'$-set.

We partition each $\xi\in \Xi'$ into a disjoint union $\xi=
\bigcup_{m\in\mathbb{Z}} \xi_m$ of $P$-towers $\xi_m$ such that the base of
$\xi_m$ is uncountable. To define $S$ on $\xi$, we apply the method used in
the proof of  \ref{smooth}. For fixed $\xi_m,\ m\in \mathbb Z,$, we set
$S =P$ everywhere except the top of $\xi_m$ and  set $S= R_m$ on the top
where $R_m$ is a Borel isomorphism mapping the top of $\xi_m$ onto the base
of $\xi_{m+1}$. Then $S$ is defined everywhere on $X$ and is aperiodic. Note
that every $\xi\in \Xi'$ is $S$-invariant. Since all towers $\xi_m$ are of the
same height, we have that $S$ maps $\xi$-atoms onto themselves. It follows
from this observation that $SF_i = PF_i,\ i=1,...,n$, that is $S\in W_0$.
\hfill$\square$

\begin{corollary}\label{xx} $\overline{\ap}^{p_0}=\auo$.
\end{corollary}
{\it Proof}. This result is an easy consequence of  \ref{dense}.
\hfill$\square$

%%%%%%%%   The Rokhlin property %%%%%%%%%%%%%%%%%%%%%%%%%%%%%

\coun The next statement proves a Borel version of the Rokhlin property for
$(\auo,p_0)$. Note that this property was considered in the settings of
measurable and Cantor dynamics in \cite{GK,GW,BDM,R}.

\begin{corollary}{\rm (the Rokhlin property)} The action of $\auo$ on itself
by conjugation is transitive with respect to the topology $p_0$. Moreover,
$\auo = \overline{\{T^{-1}ST\,:\,T\in\auo\}}^{p_0}$   holds for any
$S\in \mathcal{S}\cap\ap$.
\end{corollary}
{\it Proof}. To prove this result it suffices to use  \ref{xx}
together with the simple fact that any two aperiodic smooth automorphisms
are conjugate.

\coun The famous Rokhlin lemma on approximation aperiodic automorphisms in
the uniform topology was proved in the context of Borel dynamics in
\cite{BDK1,KM,Nad}. We formulate here this statement in the following form.

%%%%%%%%%%%5 Rokhlin Lemma %%%%%%%%%%%%%%%%%%%%%%%%%%%%%%%%55

\begin{proposition}\label{Rokhlin}  Let $T$ be an aperiodic Borel
automorphism of a standard Borel space $(X,\mathcal{B})$ and let
$\mu_1,\ldots,\mu_k\in\mathcal{M}_1(X),$ $\e>0$, and $n,m\geq 2$. Then there
exists a Borel partition of $X$ into $T$-towers $\Xi=\{\xi_k : k\in\N\}$
such that the following properties hold:

(i) The height $h_\xi$ of each $T$-tower $\xi\in\Xi$ is greater
than $n+m$.

(ii)
$\mu_j\left(\bigcup\limits_{\xi\in\Xi}\left(\bigcup\limits_{i=0}^{n-1}
T^iB_\xi \cup\bigcup\limits_{i=1}^mT^{h_\xi-i} B_\xi\right)\right)<\e$,
$j=1,\ldots,k$, where $B_\xi$ is the base of $\xi\in\Xi$.
\end{proposition}
{\it Proof.} The proof can be deduced from \cite[Theorem 2.5]{BDK1}.
\hfill$\square$

\begin{theorem}
Let $S\in\ap$. Then for any $R\in\ap$ and any $\tau$-neighborhood
$U(R)=U(R;\mu_1,\ldots,\mu_p;\e)$, there exists  $T\in\au$
such that $T^{-1}ST\in U(R)$. In other words,
$\overline{\{T^{-1}ST\,:\,T\in\au\}}^\tau=\ap$.
\end{theorem}
{\it Proof.} Apply  \ref{Rokhlin} for $R\in\ap$, $\mu_1,\ldots,
\mu_p\in\mathcal{M}_1(X)$, $\e/2>0$, and $n=m=2$. We obtain a partition of
$X$ into $R$-towers $\Xi=\{\xi_k : k\in\N\}$ satisfying (i), (ii). Choose a
sufficiently large $K\in \mathbb N$ such that $\mu_j(\bigcup\limits_{k>
K}C_k)<\e/2$, $j=1,\ldots,p$, where $C_k$ is the set supporting the tower
$\xi_k$. Therefore, we have that
$$\mu_j\left(\bigcup\limits_{k=1}^{K}\bigcup\limits_{i=2}^{h_k-3} T^iB_k
\right)>1-\e
$$ for all  $j=1,\ldots,p$, where $B_k$ is the base of $\xi_k$ and
$h_k$ is its height.

Since $S$ is aperiodic, we can find $K$ disjoint $S$-towers
$\Lambda=\{\lambda_l : k=1,\ldots, K\}$ such that the height of $\lambda_k$
is $h_k$ and $\lambda_k$ has the same cardinality as $\xi_k$ for all
$k=1,\ldots,K$. Denote by $Z_k$ the base of $\lambda_k$ and let $D_k=
\bigcup_{i=0}^{h_k-1}T^iZ_k$ be the support of $\lambda_k$. Let $Q_k$ be a
Borel isomorphism which maps $B_k$ onto $Z_k$, $k=1,\ldots,K$, and let $Q$
be a Borel isomorphism which sends $X- (C_1\cup\ldots\cup C_K)$ onto $X-
(D_1\cup\ldots\cup D_K)$. Define the automorphism $T$ as follows:
$$Tx=\left\{
\begin{array}{ll}
S^iQ_kR^{-i}, & \mbox {if}\ x\in R^iB_k,\ 0\leq i\leq h_k-1,\
1\leq k\leq K\\
\\
Qx, & \mbox {if}\ x\notin C_1\cup\ldots\cup C_K
\end{array}
\right.
$$
Then $T$ is defined everywhere on the set $X$. It
is not hard to see that $$
\{x\in X : Rx= T^{-1}STx\ \mbox {and}\ R^{-1}x= T^{-1}S^{-1}Tx\}\supset
\bigcup\limits_{k=1}^K\bigcup\limits_{i=2}^{h_k-3}R^iB_k.
$$
Hence $\mu_j(E(R,T^{-1}ST)) < \e$, $j=1,\ldots,p$, and therefore
$T^{-1}ST\in U(R)$. \hfill$\square$

%%%%%%%%%%%%%%%%%%%%%%%%%%%%%%%%%%%%%%%%%%%%%%%%%%%%%%%%%%

%%   PATH-CONNECTEDNESS

\section{Path-connectedness of $(\auo, \tau_0)$}

\coun In this section, we  prove that $\auo$ is path-connected in the
topology $\tau_0$. We first show that the group $\au$ does not possess this
property.

\begin{proposition}\label{non-path-connectness}
 The topological group $(\au,\tau)$ is not path-connected.
\end{proposition}
{\it Proof.} Let $P$ be an arbitrary involution in $\au$, that is $Px\neq x$
and $P^2x=x$ for all $x\in X$. We will show that $P$ cannot be connected
with the identity $\mathbb{I}$ by a continuous path, i.e. there exists no
continuous map $f : [0,1]\to \au$ such that $f(0)= \mathbb I, f(1) = P$.
Assume that the converse is true and let $f$ be such a path. Choose
$x_0,y_0$ in $X$ such that $Px_0 = y_0, Py_0 = x_0$. Consider the
$\tau$-neighborhood $U(P) =U(P;\delta_{x_0},\delta_{y_0};1/2)$ of $P$.
Notice that $U(P)$ contains only those automorphisms from $\au$ which map
$x_0$ to $y_0$ and $y_0$ to $x_0$. Since, by assumption, $f$ is continuous,
there exists $t^*\in (0,1)$ such that $f((t^*,1])\subset U$. Set
$$
t_0^*=\inf\{t\in [0,1]\,:\, f(s)\in U,\;t\leq s\leq 1\}.
$$
Clearly, $0\leq t_0^*<1$. Consider now the neighborhood $U(f(t_0^*)) =
U(f(t_0^*);\delta_{x_0},\delta_{y_0};1/2)$ of $f(t_0^*)$. If $t_0^*> 0$,
then there exist $\alpha$ and $\beta$ such that $\alpha< t_0^*<\beta$ and
$f([\alpha,\beta])\subset U(f(t_0^*))$. We obtain that $f(\beta)\in
U(f(t_0^*))\cap U(P)$ and therefore, $f(t_0^*)x_0=f(\beta)x_0=Px_0= y_0$ and
$f(t_0^*)y_0=f(\beta)y_0=Py_0=x_0$. A similar relation holds for
$f(\alpha)$. Thus, $f(t_0^*)\in U(P)$ and therefore $t_0^*$ cannot be the
infimum. Hence $t_0^*=0$ and $\mathbb I= f(0)\in U(P)$, which is a
contradiction. \hfill$\square$

\begin{remark}\label{sect}  {\rm Let $T\in \au$ and let $A\in \B$ be a
complete $T$-section such that every point from $A$ is recurrent. If a Borel
set $B$ contains $A$, then $B$ is also a complete $T$-section which consists
of recurrent points.}
\end{remark}

\begin{theorem} The topological group  $(\auo,\tau_0)$ is path-connected.
\end{theorem}
{\it Proof}. We first prove separately that every periodic automorphism $P$
and every aperiodic automorphism $T$ can be connected with the identity by a
continuous path (see \ref{per} and \ref{aper} respectively). By
\ref{technicals}, these two results will give the proof for any
automorphism. Recall that, by \ref{coun}, it is sufficient to deal with
continuous measures only.

%%%   Connecting a periodic automorphism  with the identity

\begin{lemma}\label{per} Let $P\in\auo$ be a periodic automorphism. Then
there exists a continuous map $f:[0,1]\rightarrow (\auo,\tau_0)$ such that
$f(0)=\mathbb{I}$ and $f(1)=P$.
\end{lemma}
{\it Proof.} By \ref{technicals}, we have the decomposition of
$X=\bigcup_{k\geq 1}X_k$ where $X_k=\bigcup_{i=0}^{k-1}P^iB_k$ is a
$P$-tower. Without loss of generality, we can assume that the $B_k$'s are
uncountable Borel sets and therefore they all are isomorphic to the unit
interval $(0,1)$. Let $\psi_k : (0,1)\rightarrow B_k,\ k\in \mathbb N$ be a
Borel isomorphism. For each $B_k$, define the map $\Psi_k : [0,1]\rightarrow
\mathcal{B}\upharpoonright B_k$ by
$$
\Psi_k(t)=\left\{
\begin{array}{ll}
\emptyset, & t=0\\
\psi_k((0,t)), & t\in(0,1)\\
B_k, & t=1.
\end{array}
\right.
$$

\noindent {\it Claim 1}. The function $t\mapsto \Psi_k(t),\ k\in \mathbb
N,$ is continuous on [0,1] in the sense that for any non-atomic
$\mu\in\mathcal{M}_1(X)$,
$$
\lim\limits_{t\to t_0 }\mu(\Psi_k(t)\ \triangle\ \Psi_k(t_0))=0,\ \ t_0\in [0,1].
$$
The proof is straightforward.
\medskip

Define now the path $f:[0,1]\rightarrow \auo$ as follows:
$$
f(t)x=\left\{
\begin{array}{ll}
Px, & \mbox {if} \ x\in \bigcup\limits_{i=0}^{k-1}P^i\Psi_k(t)
\ \mbox {for some} \ k\\
x, & \mbox{otherwise}.
\end{array}
\right.
$$
It is clear that $f(0)=\mathbb{I}$ and $f(1)=P$ and we need to show only that
$f(t)$ is continuous. To do this, fix $t_0\in [0,1]$ and consider the map
$\Theta: t \mapsto \mu(E_0(f(t),f(t_0)))$ on [0,1] where
$\mu\in\mathcal{M}_1(X)$ is non-atomic and $E_0$ is defined in \ref{E0}.

Given $\e>0$, choose $K>0$ such that $\mu(\bigcup_{k>K}X_k)<\e/2$. Therefore
$\mu(E_0(f(t),f(t_0)))\leq \mu\left(E_0(f(t),f(t_0))\cap
\bigcup_{k=1}^K X_k\right)+\e/2$. We see that
$$
\mu\left(E_0(f(t),f(t_0))\cap \bigcup_{k=1}^KX_k\right)=
\sum\limits_{k=1}^K\sum\limits_{i=0}^{k-1}\mu\left(P^i(\Psi_k(t)\triangle
\Psi_k(t_0) )\right).
$$
The fact that $\Theta$ is continuous follows from Claim 1.

If now $U_0(f(t_0);\mu_1,\ldots,\mu_n;\e)$ is a $\tau_0$-neighborhood of
$f(t_0)$, then we apply the proved result to each measure $\mu_i$. The lemma
is proved. \hfill$\square$

%%%%%%%%   connecting an aperiodic automorphism with the identity %%%%%%

\begin{lemma}\label{aper} Let $T\in \auo$ be an arbitrary aperiodic
automorphism. Then there exists a continuous map $P: [0,1] \rightarrow
(\auo,\tau_0)$ such that $P(0)=\mathbb{I}$ and $P(1)=T$. Moreover, for all
$t\neq 1$, the automorphism $P(t)$ is periodic.
\end{lemma}
{\it Proof.} By  \ref{markers}, choose a vanishing sequence of markers
$\{A_n\}_{n=0}^\infty$ with $A_0=X$. Without loss of generality, we can
assume that the set $F_n:=A_{n}\setminus A_{n+1}$ is uncountable for all
$n$. Take a sequence of real numbers $\{t_n\}$ such that $0=t_0<t_1<t_2<
\ldots <1$ and $\lim\limits_{n\to\infty} t_n=1$. Let $\psi_n :[t_n,
t_{n+1})\rightarrow F_n,\ n\in \mathbb N,$ be a Borel isomorphism. Define
the function $\Phi : [0,1] \rightarrow \mathcal{B}$ as follows:
$$
\Phi(t)=\left\{
\begin{array}{lll}
A_n, & \mbox{if} \ t=t_n,\ n\geq 0\\
A_n-\psi_n([t_n,t)), & \mbox{if} \ t\in (t_n,t_{n+1}),\ n\geq 0\\
\emptyset, & \mbox{if}\ t=1.
\end{array}\right.
$$
Observe that for each
$t\in[0,1)$ there exists $n= n(t)\in\N$ such that $\Phi(t)\supset A_n$. By
 \ref{sect}, we get that $\Phi(t)$ is a $T$-section which consists of
recurrent points. We also notice that $\lim\limits_{t\to
s}\mu(\Phi(t)\triangle \Phi(s))=0$ for any non-atomic
$\mu\in\mathcal{M}_1(X)$ and $s\in[0,1]$.

Now, we apply the method of the proof of the Rokhlin lemma \cite{BDK1} to
produce a continuous family $\{P(t)\}$ of periodic automorphisms which
approximates $T$.

Fix $t\in[0,1]$. By \ref{Towers_Construction}, define the set
$\Phi_n(t)=\{x\in\Phi(t) : T^nx\in \Phi(t), T^ix\notin \Phi(t) \mbox{ for
}1\leq i\leq n-1 \},\ n\in \mathbb N$, where $\Phi_n(1)=\emptyset$ and
$\Phi_n(0)=X$. Clearly, $\Phi_n(t)$ may be empty for some $n$. By
construction, the entire space $X$ is partitioned into $T$-towers with bases
$\Phi_n(t)$.

Define
$$
P(t)=\left\{
\begin{array}{ll}
T^{-n+1}x, & \mbox{whenever }x\in T^{n-1}\Phi_n(t) \mbox{ for some
}n\in\N\\
Tx, & \mbox{otherwise}.
\end{array}\right.
$$
Notice that $P(0)=\mathbb{I}$, $P(1)=T$, and  $P(t)$ is periodic if $t\neq 1$.

\par\medskip\noindent
It remains to prove that the map $P: t\mapsto P(t)$ sending $[0,1]$ to
$\auo$ is continuous. Note that for every $n\in\N$ the family
$\{\Phi_n(t)\}$ is continuous, i.e. for a non-atomic measure $\mu\in\M$,
$$
\lim\limits_{t\to s} \mu(\Phi_n(t)\triangle \Phi_n(s))=0.
$$
Indeed, this fact follows from continuity of $\Phi(t)$ and from the relation
$$
\Phi_n(t)=(\Phi(t)\cap T^{-n}\Phi(t))\setminus
\bigcup\limits_{i=1}^{n-1}T^{-i}\Phi(t).
$$

Next, show that for $t',t'' \in [0,1]$ we have
$$
E_0(P(t'),P(t''))=\{x\in X : P(t')x\neq P(t'')x\}\subset
\bigcup\limits_{n=1}^\infty\bigcup\limits_{k=0}^{n-1}T^k(\Phi_n(t')\triangle
\Phi_n(t'')).
$$
Let $t'<t''$ for definiteness. Suppose that $x\in E_0(P(t'),P(t''))$. Then
$x$ belongs to a $T$-tower constructed over $\Phi(t')$, that is $x\in
T^l\Phi_n(t')$ for some $n\in\N$ and $0\leq k\leq n-1$. If $x\in
T^k\Phi_n(t'')$, then by construction of $P(t')$ and $P(t'')$, we have that
$P(t')x=P(t'')x$. Therefore, $x\in T^k(\Phi_n(t')\setminus \Phi_n(t''))$.
\medskip

Fix $s\in[0,1]$. Consider a neighborhood $U(P(s))=U(P(s);
\mu_1,\ldots,\mu_m; \e)$ of $P(s)$ where measures $\mu_1,\ldots,\mu_m$ are
continuous. By definition of $P(s)$, we have that
$$
X=\bigcup\limits_{n\geq
1}\bigcup\limits_{k=0}^{n-1}P^k(s)\Phi_n(s)=
\bigcup\limits_{n\geq 1}\bigcup\limits_{k=0}^{n-1}T^k\Phi_n(s),
$$
and these unions are disjoint. Take $N\geq 1$ such that
$$
\mu_i\left(\bigcup\limits_{n\geq
N}\bigcup\limits_{k=0}^{n-1}T^k\Phi_n(s)\right)<\e/8,\qquad (*)
$$
for $1\leq i\leq m$. By continuity of $\{\Phi_n(t)\}$, we can find a
neighborhood $O(s)\subset [0,1]$ such that for any $t\in O(s)$ one has
$$
\mu_i\left(T^k\left(\Phi_n(t)\triangle \Phi_n(s)\right)\right)<
\frac{\e}{8N^2}\qquad (**)
$$
for $1\leq i\leq m$, $1\leq n\leq N-1$, and $0\leq k\leq n-1$.

It follows from (*) and (**) that
$$
\mu_i\left(\bigcup\limits_{n=1}^{N-1}\bigcup\limits_{k=0}^{n-1}
T^k\Phi_n(t)\right)>1-\e/4,
$$
hence $\mu_i(\bigcup\limits_{n\geq
N}\bigcup\limits_{k=0}^{n-1}T^k\Phi_n(t))<\e/4$ when $i=1,\ldots,m$ and
$t\in O(s)$.

Therefore, we have
$$\begin{array}{lll}
\mu_i\left(\bigcup\limits_{n\geq
N}\bigcup\limits_{k=0}^{n-1}T^k(\Phi_n(t)\triangle
\Phi_n(s))\right)\\
\leq \mu_i\left(\bigcup\limits_{n\geq
N}\bigcup\limits_{k=0}^{n-1}T^k
\Phi_n(s)\right)+\mu_i\left(\bigcup\limits_{n\geq
N}\bigcup\limits_{k=0}^{n-1}T^k\Phi_n(t)\right)\\<\e/2,
\end{array}
$$
for all $t\in O(s)$.
Thus, for all $t\in O(s)$, we obtain that
$$
\begin{array}{lll}
\mu_i(E_0(P(t),P(s)))\leq
\mu_i\left(\bigcup\limits_{n=1}^\infty\bigcup\limits_{k=0}^{n-1}
T^k(\Phi_n(t)\triangle \Phi_n(s))\right)\\
\\
\leq\mu_i\left(\bigcup\limits_{n=1}^{N-1}\bigcup\limits_{k=0}^{n-1}
T^k(\Phi_n(t)\triangle \Phi_n(s))\right) +
\mu_i\left(\bigcup\limits_{n=N}^\infty\bigcup\limits_{k=0}^{n-1}
T^k(\Phi_n(t)\triangle \Phi_n(s))\right)\\
\\
\displaystyle\leq
\sum\limits_{n=1}^{N-1}\sum\limits_{k=0}^{n-1}\frac{\e}{8N^2}+
\frac{\e}{2}\\
\\
<\e.
\end{array}$$
This means that $\{P(t): t\in O(s)\}\subset U(P(s))$ and the proof is
completed. \hfill$\square$
\\

\coun Let $T\in \au$ and let $Orb_T(x)$ denote the $T$-orbit of $x\in X$.
Recall the definition of the {\it full group} $ [T] $ generated by $T\in
\au$:
$$
[T] = \{ \gamma \in \au\ \vert\ \gamma x \in Orb_T(x),\ \forall x\in X\}.
$$
Then every $\gamma \in [T]$ defines a Borel function $m_\gamma : X \to
\mathbb Z$ such that $\gamma x = T^{m_\gamma (x)}x,\ x\in X $. It follows
easily from  \ref{top} that $[T]$ is $\tau$-closed in $\au$.

Note that if $T\sim S$, then $T^n \sim S^n,\ \forall n\in {\mathbb Z}$.
Therefore $Orb_T(x) = Orb_S(x)$ everywhere except a countable set. This
means that one can extend the definition of full group to automorphisms from
$\auo$.

\begin{corollary} The full group $[T]$ of any $T\in \auo$ is path-connected.
\end{corollary}
{\it Proof}. The proofs of  \ref{per} and \ref{aper} show that the
constructed paths $f(t)$ and $P(t)$ connecting the identity with $T$ belong
to the full group $[T]$. \hfill$\square$
\\

{\it Acknowledgement}. We are thankful to A.H.~Dooley, A.~Kechris,
J.~Kwiatkowski, and B.~Miller for numerous discussions of these results. The
first-named author thanks the University of New South Wales for the warm
hospitality and the Australian Research Council for its support.

%%%%%%%%%%%%%%%%%%%%%%%%%%%%%%%%%%%%%%%%%%%%%%%%%%%%%%%%%%

%%%%%%%%    REFERENCES
%\newpage

\vskip1cm \noindent {\small {\it S.~Bezuglyi and K.~Medynets\\
Institute for Low Temperature Physics\\
47 Lenin ave., 61103 Kharkov\\
Ukraine\\

\noindent bezuglyi@ilt.kharkov.ua\\
medynets@ilt.kharkov.ua} }
\end{document}